\documentclass[11pt]{article}
\usepackage{graphicx}
\usepackage{amsmath}
\usepackage{amsfonts}
\usepackage{amssymb}
\newtheorem{theorem}{Theorem}

\newtheorem{definition}[theorem]{Definition}

\newtheorem{lemma}[theorem]{Lemma}

\newtheorem{proposition}[theorem]{Proposition}

\newenvironment{proof}[1][Proof]{\textbf{#1.} }{\ \rule{0.5em}{0.5em}}

\begin{document}

\title{2D models of statistical physics with continuous symmetry: the case of
singular interactions.}
\author{D. Ioffe\\Technion, Haifa\\\textit{ieioffe@ie.technion.ac.il}
\and S. Shlosman\\CPT, Luminy, Marseille \\\textit{shlosman@cpt.univ-mrs.fr}
\and Y. Velenik\\UMR-CNRS 6632, CMI, Marseille\\\textit{velenik@cmi.univ-mrs.fr}}
\maketitle
\begin{abstract}
We show the absence of continuous symmetry breaking in 2D lattice systems
without any smoothness assumptions on the interaction. We treat certain cases
of interactions with integrable singularities. We also present cases of
singular interactions with continuous symmetry, when the symmetry is broken in
the thermodynamic limit.

\end{abstract}

\section{Introduction and results}

\subsection{The invariance problem: an overview}

\noindent In this paper we are studying the two-dimensional lattice models of
statistical mechanics, which are defined by a $G$-invariant interaction, where
$G$ is some compact connected Lie group. We shall investigate both the cases
of finite and infinite range interactions. The general class of finite-range
models to be considered is given by the following Hamiltonians:
\begin{equation}
\mathcal{H}\left(  \mathbf{\phi}\right)  =\sum_{x\in\mathbb{Z}^{2}}U_{\Lambda
}\left(  \phi_{\cdot+x}|_{\Lambda}\right)  . \label{01}
\end{equation}
Here $\mathbf{\phi}=\left\{  \phi_{y},y\in\mathbb{Z}^{2}\right\}  $ is the
field, taking values in some compact topological space $S,$ $\Lambda$ is a
fixed finite subset of $\mathbb{Z}^{2}$ and the translation-invariant
interaction\textbf{\ }$\mathcal{U}=\left\{  U_{\Lambda+x}\left(  \cdot\right)
,x\in\mathbb{Z}^{2}\right\}  $\textbf{\ }is specified by a real function
$U_{\Lambda}$\ on $S^{\Lambda}.$ We suppose that a continuous action of a
compact connected Lie group $G$ on $S$ is given, $a:G\times S\rightarrow S,$
and for $\phi\in S,$ $g\in G$ we introduce the notation $g\phi=a\left(
g,\phi\right)  .$ This action defines the action of $G$ on $S^{k}$ for every
$k$ by $g\left(  \phi_{1},...,\phi_{k}\right)  =\left(  g\phi_{1}
,...,g\phi_{k}\right)  ,$ and the main assumption is that the function
$U_{\Lambda}$ is invariant under this action on $S^{\Lambda}$: for every $g\in
G$
\begin{equation}
U_{\Lambda}\left(  g\left(  \phi_{1},...,\phi_{\left|  \Lambda\right|
}\right)  \right)  =U_{\Lambda}\left(  \phi_{1},...,\phi_{\left|
\Lambda\right|  }\right)  . \label{04}
\end{equation}
Of course, we suppose that the free measure $d\phi$ on $S$ is $G$-invariant as well.

The best known examples of such models are XY model (or plane rotator model)
and XYZ model (or classical Heisenberg model). For the XY model $S=\mathbb{S}
^{1}\subset\mathbb{R}^{2}$ is the unit circle, and for the XYZ model
$S=\mathbb{S}^{2}\subset\mathbb{R}^{3}$ is the unit sphere. The Hamiltonians
are
\begin{equation}
\mathcal{H}\left(  \mathbf{\phi}\right)  =-\sum_{\substack{x,y\in
\mathbb{Z}^{2} \\\left|  x-y\right|  =1}}J\left(  \phi_{x}\circ\phi
_{y}\right)  , \label{02}
\end{equation}
where $J$ is a positive constant, $\left(  ..\circ..\right)  $ stands for the
scalar product, and the free measures are just the Lebesgue measures on the spheres.

The first rigorous result in this field is the well-known Mermin-Wagner
theorem, which for the models (\ref{02}) states the absence of spontaneous
magnetization. Then in \cite{DS1} a stronger result was proven, stating that
under some \textbf{smoothness conditions} (see (\ref{03}) below) on the
function $U_{\Lambda}$ every Gibbs state of the model defined by the
Hamiltonian (\ref{01}) is $G$-invariant under the natural action of $G$ on
$S^{\mathbb{Z}^{2}}.$ Later it was proven in \cite{MS} that for the model
(\ref{02}) the correlations decay at least as a power law. In \cite{S78} the
same power law was obtained for the general model (\ref{01}), again under the
smoothness condition (\ref{03}). This result was reproved later in \cite{N}
for the case $G=SO(n),$ by means of the complex translations method of McBryan
and Spencer. Another proof of $G$-invariance of the Gibbs states of models of
type (\ref{01}) was found in~\cite{P, FP}; with this technique it was possible
to prove the result for long range interactions decaying as slowly as $r^{-4}$
(in fact better, see the remark after Theorem~\ref{thm2}). This result is
optimal since it is known~\cite{KP} that in the low temperature XY model with
interactions decaying as $r^{-4+\alpha}$, $\alpha>0$, there is spontaneous
symmetry breaking. On the other hand, this technique seems unable to yield the
algebraic decay of correlations. In the case of $SO(N)$-symmetric models, the
technique of~\cite{MS} can be extended to cover such long-range interactions,
see~\cite{MMR}.

An alternative approach to these problems is via Bogoliubov inequalities,
see~\cite{M, KLS, BPK, I}; it also permits to prove absence of continuous
symmetry breaking for long-range interactions. As the technique of~\cite{P}
however, they seem unable to yield the correct decay of correlations.

The smoothness condition on the interaction, which was crucial for all the
results mentioned above, is the following. Let $M\subset\Lambda$ be a subset,
and $\phi_{\Lambda}=\left(  \phi_{M}\cup\phi_{\Lambda\backslash M}\right)  $
be an arbitrary configuration. Then one requires that for every choice of the
subset $M$ and the configuration $\phi_{\Lambda}$ the functions\footnote{In
fact, as noted in~\cite{P}, the proofs really only use the fact that these
functions are $\mathcal{C}^{1}$, with a first derivative satisfying a
Lipschitz condition.}
\begin{equation}
V_{\phi_{\Lambda},M}\left(  g\right)  =U_{\Lambda}\left(  \phi_{M}\cup
g\phi_{\Lambda\backslash M}\right)  \text{ are }\mathcal{C}^{2}\text{
functions on }G. \label{03}
\end{equation}
Moreover, the second derivatives of the functions $V_{\phi_{\Lambda},M}\left(
\cdot\right)  ,$ taken along any tangent direction in $G,$ have to be bounded
from above, uniformly in $\phi_{\Lambda}$ and $M$. In the next section we are
explaining how that condition can be used in the proof of the $G$-invariance.

\subsection{No breaking of continuous symmetry for singular interactions}

\textbf{The Main Result} \textbf{of the present paper} is that the smoothness
property (\ref{03}) is in fact not necessary for the $G$-invariance, and it is
implied by the mere continuity of the functions $V_{\phi_{\Lambda},M}\left(
g\right)  $ and the invariance (\ref{04}). Moreover, even the continuity is
not necessary, and a certain integrability condition on $V_{\phi_{\Lambda}
,M}\left(  g\right)  $ is enough (see relation (\ref{51}) below). For example,
for $G$ being a circle, $\mathbb{S}^{1},$ with $U$ nearest neighbour
interaction, $U\left(  \phi_{1},\phi_{2}\right)  =U\left(  \phi_{1}-\phi
_{2}\right)  ,$ the singularity $U\left(  \phi\right)  \sim\ln\left|
\phi\right|  $ at $\phi=0$ does not destroy $\mathbb{S}^{1}$-invariance of the
corresponding Gibbs measures. However, if the interaction $U$ is ``even more
singular'', then the $G$-invariance can be destroyed, as the Theorem \ref{Aiz}
below shows.

For simplicity we will consider the case when both the space $S$ and the group
$G$ will be a circle, $\mathbb{S}^{1}.$ The general case follows easily from
this special one, see \cite{DS1}, since for every element $g\in G$ there is a
compact commutative subgroup (torus) $T\subset G,$ such that $g\in T.$ We also
suppose that the interaction\textbf{\ }$\mathcal{U}$ is a nearest neighbour
translation invariant interaction, given by a symmetric function $U$ of two
variables: $U\left(  \phi_{1},\phi_{2}\right)  =U\left(  \phi_{2},\phi
_{1}\right)  .$ The generalization to a finite range interaction is
straightforward. The $\mathbb{S}^{1}$-invariance of $\mathcal{U}$ means that
$U\left(  \phi_{1},\phi_{2}\right)  =U\left(  \phi_{1}+\psi,\phi_{2}
+\psi\right)  $ for every $\psi\in\mathbb{S}^{1},$ so in fact we can say that
$U$ is a function of one variable, $U\left(  \phi_{1},\phi_{2}\right)
=U\left(  \phi_{1}-\phi_{2}\right)  ,$ with $U\left(  \phi\right)  =U\left(
-\phi\right)  .$

To formulate our main result, we will introduce the notation $\mu^{U_{\Lambda
}}$ for a Gibbs measure on $S^{\mathbb{Z}^{2}}$ corresponding to the formal
Hamiltonian (\ref{01}) and the free measure $d\phi$ on $S,$ which is supposed
to be $G$-invariant. We will denote the integration operation with respect to
$\mu^{U_{\Lambda}}$ by $\left\langle \cdot\right\rangle ^{U_{\Lambda}}.$ By
$\mu_{x,y}^{U_{\Lambda}}\left(  d\phi_{x},d\phi_{y}\right)  $ we denote the
restriction of the measure $\mu^{U_{\Lambda}}$ to the product $S\times S$ of
the state spaces of the field variables $\phi_{x}$ and $\phi_{y}.$

We will prove the following

\begin{theorem}
Suppose that:

\begin{itemize}
\item  the finite range interaction function $U_{\Lambda}$ is continuous,
bounded on $S^{\left|  \Lambda\right|  },$ and satisfies the $G$-invariance
property (\ref{04}),

\item  the free measure $d\phi$ on $S$ is $G$-invariant.

Then the measure $\mu^{U_{\Lambda}}$ is $G$-invariant: for every $g\in G,$
every finite $V\subset\mathbb{Z}^{2}$ and every $\mu^{U_{\Lambda}}$-integrable
function $f$ on $S^{V}$
\begin{equation}
\left\langle f\left(  g\cdot\right)  \right\rangle ^{U_{\Lambda}}=\left\langle
f\left(  \cdot\right)  \right\rangle ^{U_{\Lambda}}. \label{06}
\end{equation}
Moreover, it has the following correlation decay: for every $A,B\subset S$ the
conditional distributions of the measure $\mu_{x,y}^{U_{\Lambda}}$ satisfy for
every $g\in G$ the estimate
\begin{equation}
\left|  \frac{\mu_{x,y}^{U_{\Lambda}}\left(  \phi_{x}\in gA{\Huge \,}\Bigm
|\,\phi_{y}\in B\right)  }{\mu_{x,y}^{U_{\Lambda}}\left(  \phi_{x}\in
A{\Huge \,}\Bigm|\,\phi_{y}\in B\right)  }-1\right|  \leq C\left(  U_{\Lambda
}\right)  \left|  x-y\right|  ^{-c\left(  U_{\Lambda}\right)  }, \label{05}
\end{equation}
with $C\left(  U_{\Lambda}\right)  <\infty,$ $c\left(  U_{\Lambda}\right)
>0.$ In case when the space $S$ is a homogeneous space of the group $G$ (e.g.
$S=G$), the measure $\mu_{x,y}^{U_{\Lambda}}\left(  d\phi_{x},d\phi
_{y}\right)  $\ can be written as a convex sum of two probability measures:
\[
\mu_{x,y}^{U_{\Lambda}}\left(  d\phi_{x},d\phi_{y}\right)  =c_{xy}\hat{\mu
}_{x,y}^{U_{\Lambda}}\left(  d\phi_{x},d\phi_{y}\right)  +\left(
1-c_{xy}\right)  \tilde{\mu}_{x,y}^{U_{\Lambda}}\left(  d\phi_{x},d\phi
_{y}\right)  .
\]
The measure $\hat{\mu}_{x,y}^{U_{\Lambda}}\left(  d\phi_{x},d\phi_{y}\right)
$\ can be singular, but the number $c_{xy}$\ is very small: $0\leq c_{xy}
\leq\exp\left\{  -\sqrt{\left|  x-y\right|  }\right\}  ,$\ while the measure
$\tilde{\mu}_{x,y}^{U_{\Lambda}}\left(  d\phi_{x},d\phi_{y}\right)  $\ has a
density $p_{x,y}^{U_{\Lambda}}\left(  \phi_{x},\phi_{y}\right)  $\ with
respect to the measure\textbf{\ }$d\phi_{x}d\phi_{y},$ which for every
conditioning $\phi_{y}=\psi$ satisfies the estimate
\[
\left|  p_{x,y}^{U_{\Lambda}}\left(  \phi_{x}\Bigm|\phi_{y}=\psi\right)
-1\right|  \leq C\left(  U_{\Lambda}\right)  \left|  x-y\right|  ^{-c\left(
U_{\Lambda}\right)  },
\]
with $c\left(  U_{\Lambda}\right)  >0.$ In particular, for the case
$G=SO\left(  n\right)  ,$ $S=\mathbb{S}^{n-1}\subset\mathbb{R}^{n}\ $with
$n\geq2$ we have
\begin{equation}
0\leq\left\langle \left(  \phi_{x}\circ\phi_{y}\right)  \right\rangle
^{U_{\Lambda}}\leq C\left(  U_{\Lambda}\right)  \left|  x-y\right|
^{-c\left(  U_{\Lambda}\right)  }. \label{07}
\end{equation}
\end{itemize}
\end{theorem}

We remind the reader that the homogeneous space is a manifold of the classes
of conjugacy of a compact subgroup $H\subset G.$

The $G$-invariance (\ref{06}) does not imply the uniqueness of the Gibbs state
with the interaction $U_{\Lambda}.$ The reason is that the interaction
$U_{\Lambda}$ may possess an additional discrete symmetry, which may be
broken. An example is constructed in \cite{S80}.

The estimate (\ref{07}) cannot be improved in general. Indeed, Fr\"{o}hlich
and Spencer have obtained the power law decay of the pair correlations in XY
model (\ref{02}) for large values of the coupling constant $J,$ see \cite{FS}.
On the other hand, for XYZ model it is expected that the pair correlations
decay exponentially for all values of $J.$

\subsection{Infinite range case}

The preceding theorem is restricted to finite-range interactions. Let us now
turn to the long-range case. The formal Hamiltonian is supposed to be of the
form
\begin{equation}
\mathcal{H}\left(  \mathbf{\phi}\right)  =\sum_{x,y}J_{x-y}\,U(\phi_{x}
,\phi_{y})\,. \label{70}
\end{equation}
More general Hamiltonians (e.g., without separating the spatial and spin part
of the interaction, or with more than 2-body interactions) could also be
treated along the lines of the approach we develop here, but for the sake of
simplicity we shall restrict ourselves to the case of (\ref{70}). Since the
coupling constants $\{J_{\cdot}\}$ have to satisfy the summability condition,
we can make an additional normalization assumption
\begin{equation}
\sum_{x\neq0}|J_{x}|=1. \label{J_normalized}
\end{equation}
Let $X_{\cdot}$ be the random walk on $\mathbb{Z}^{2}$ with transition
probabilities from $x$ to $y$ given by $|J_{x-y}|.$

We then have the following

\begin{theorem}
\label{thm2}Suppose that

\begin{itemize}
\item  The random walk $X_{\cdot}$ is recurrent.

\item  The 2-body interaction function $U$ is continuous {\textbf{ }}on
$S\times S,$ and satisfies the invariance property (\ref{04}).

\item  The free measure $d\phi$ on $S$ is $G$-invariant.
\end{itemize}

Then all Gibbs states, corresponding to the Hamiltonian (\ref{70}), are $G$-invariant.
\end{theorem}

The recurrency condition is known to be optimal even in the case of smooth
$U$, in the sense that there are examples of systems for which the continuous
symmetry is broken as soon as the underlying random-walk is transient,
see~\cite{BPK} or Theorem~(20.15) in~\cite{Georgii}.

Recurrence of the underlying random-walk is not a very explicit condition.
Explicit examples have been given in~\cite{P}. Namely, it follows from the
latter that Theorem~\ref{thm2} applies if there exists $p<\infty$ such that
the coupling constants decays for large \textbf{$\Vert x\Vert_{\infty}$} at
least\textbf{ }like\textbf{\ }
\[
\Vert x\Vert_{\infty}^{-4}\log_{2}\Vert x\Vert_{\infty}\dots\log_{p}\Vert
x\Vert_{\infty}\,,
\]
where $\log_{k}x=\log\log_{k-1}x,$ and $\log_{2}x=\log\log x$. On the other
hand, it follows from~\cite{FILS} that the continuous symmetry is broken for
the low temperature XY model with coupling constants behaving, for large
$\Vert x\Vert_{\infty}$, like
\[
\Vert x\Vert_{\infty}^{-4}\log_{2}\Vert x\Vert_{\infty}\dots\left(  \log
_{p}\Vert x\Vert_{\infty}\right)  ^{1+\varepsilon}\,,
\]
for any $p<\infty$ and $\varepsilon>0$.

\subsection{Non-compact symmetry group: non-existence of 2D Gibbs states}

Finally we mention the case of connected non-compact Lie group $G$. The case
of the smooth interaction was treated in \cite{DS2}, and the corresponding
long-range result was obtained in~\cite{FP}. Technically the compact and the
non-compact cases are very similar, but the results are quite different. The
reason is that while in the compact case the Haar measure on $G$ can be
normalized to a probability measure, in the non-compact case it is not
possible. Therefore, there are no $G$-invariant probability measures on $G$
for $G$ non-compact. This is the main reason behind the result of \cite{DS2}
and~\cite{FP}: the corresponding 1D and 2D Gibbs measures do not exist.

Below we are formulating the simplest such result for the non-compact case and
singular interaction that our technique can produce. The field $\mathbf{\phi}$
will be real-valued, $G=\mathbb{R}^{1},$ and
\begin{equation}
\mathcal{H}\left(  \mathbf{\phi}\right)  =\sum_{x,y\in\mathbb{Z}^{2}}
J_{x-y}\,\bar{U}\left(  \phi_{x}-\phi_{y}\right)  , \label{81}
\end{equation}
with the function $\bar{U}$ satisfying

\begin{itemize}
\item $\bar{U}\left(  \phi\right)  =\bar{U}\left(  -\phi\right)  ,$

\item $\bar{U}\left(  \phi\right)  =U\left(  \phi\right)  -\upsilon\left(
\phi\right)  ,$ where $U$ is a $\mathcal{C}^{2}$ function with uniformly
bounded second derivative, and $\upsilon\leq\varepsilon_{0} <1 $, where
$\varepsilon_{0}$ is some technical constant, which is small,
\end{itemize}

and the coupling constants $\{ J_{\cdot}\}$ satisfy the same hypothesis as in
Theorem~\ref{thm2}\thinspace.

\begin{theorem}
There are no two-dimensional Gibbs fields, corresponding to the Hamiltonian
(\ref{81}), with interaction $\bar{U}$ and coupling constants $J_{\cdot}$ as above.
\end{theorem}

In particular, the last theorem covers the case of the (non-convex)
interactions
\[
\bar{U}\left(  \phi\right)  =\left|  \phi\right|  ^{\alpha},\;{0<\alpha\leq1},
\]
and so answers a question which was left open in the paper \cite{BLL}. In
fact, all the results of \cite{BLL} concerning the non-existence of the 2D
Gibbs fields for interactions growing at most quadratically in $\phi$ follow
from the above theorem. Notice that our techniques also allow to obtain lower
bounds with the correct behavior for the variance of the field in a finite box.

The general formulation of the above theorem and its proof will be published
in a separate paper.

\subsection{Continuous symmetry breaking in 2D}

Our results on continuous symmetry breaking are inspired by the paper of M.
Aizenman \cite{A}, where the following was proven. Consider the case when
$S=G=\mathbb{S}^{1},$ with the interaction $U\left(  \phi_{1},\phi_{2}\right)
=U\left(  \phi_{1}-\phi_{2}\right)  $ given by
\begin{equation}
U\left(  \phi\right)  =\left\{
\begin{array}
[c]{cc}
-\cos\phi & \text{ if }\left|  \phi\right|  \leq\theta,\\
+\infty & \text{ if }\left|  \phi\right|  >\theta.
\end{array}
\right.  \label{85}
\end{equation}
Then in the 2D case, the statement of \cite{A} is that the two-point pair
correlations in the state with free or periodic b.c. decay at most as a power
law, at all temperatures including infinite temperature, provided $\left|
\theta\right|  <\frac{\pi}{4}.$

It would be interesting to know whether the Gibbs states of this model with
zero b.c., i.e. $\mathbf{\phi}\equiv0,$ are $\mathbb{S}^{1}$-invariant. To the
best of our knowledge this question is open. However, one can prove the
following simple:

\begin{theorem}
\label{Aiz} Suppose that $\theta=\theta_{k}=\frac{2\pi}{k},$ for some
$k=9,10,...$ . Then there exist translation-invariant Gibbs states $\mu_{k}$,
corresponding to the interaction (\ref{85}), which are not $\mathbb{S}^{1}$-invariant.
\end{theorem}

These states possess, however, $\mathbb{Z}_{k}$-invariance, $\mathbb{Z}
_{k}\subset\mathbb{S}^{1}$.

\section{Proofs}

\subsection{Theorem 1: Smooth case.}

We begin by reminding the reader the main ideas of the proof for the case of
smooth interaction. The proof for the general case would be built upon it. We
follow \cite{DS1}, with simplifications made in \cite{Si}.

For simplicity we will consider the case when both the space $S$ and the group
$G$ will be a circle, $\mathbb{S}^{1}.$ The general case follows easily from
this special one, see \cite{DS1}, since for every element $g\in G$ there is a
compact commutative subgroup (torus) $T\subset G,$ such that $g\in T.$ We also
suppose that the interaction\textbf{\ }$\mathcal{U}$ is a nearest neighbour
translation invariant interaction, given by a symmetric function $U$ of two
variables: $U\left(  \phi_{1},\phi_{2}\right)  =U\left(  \phi_{2},\phi
_{1}\right)  .$ The generalization to a finite range interaction is
straightforward. The $\mathbb{S}^{1}$-invariance of $\mathcal{U}$ means that
$U\left(  \phi_{1},\phi_{2}\right)  =U\left(  \phi_{1}+\psi,\phi_{2}
+\psi\right)  $ for every $\psi\in\mathbb{S}^{1},$ so in fact we can say that
$U$ is a function of one variable, $U\left(  \phi_{1},\phi_{2}\right)
=U\left(  \phi_{1}-\phi_{2}\right)  ,$ with $U\left(  \phi\right)  =U\left(
-\phi\right)  .$ The smoothness we need is the following: we suppose that $U$
has the second derivative, which is bounded from above:
\begin{equation}
U^{\prime\prime}\left(  \phi\right)  \leq\bar{C}. \label{25}
\end{equation}

Let $\Lambda_{n}$ be the box $\left\{  x\in\mathbb{Z}^{2}:\left|  \left|
x\right|  \right|  _{\infty}\le n\right\}  ,$ and $\bar\phi$ be an arbitrary
boundary condition outside $\Lambda_{n}$. Let $\left\langle \cdot\right\rangle
_{n,\bar\phi}$ be the Gibbs state in $\Lambda_{n}$ corresponding to the
interaction $U$ and the boundary condition $\bar\phi.$ Let $V$ be an arbitrary
finite subset of $\mathbb{Z}^{2},$ containing the origin. Our theorem will be
proven for the interaction $U$ once we obtain the following estimate:

\begin{lemma}
For every function $f\left(  \phi\right)  =f\left(  \phi_{V}\right)  ,$ which
depends only on the configuration $\phi$ inside $V,$ we have for every
$\psi\in\mathbb{S}^{1}$
\begin{equation}
\left|  \left\langle f\left(  \phi+\psi\right)  \right\rangle _{n,\bar{\phi}
}-\left\langle f\left(  \phi\right)  \right\rangle _{n,\bar{\phi}}\right|
\leq C\left(  \bar{C},V\right)  \left|  \left|  f\right|  \right|  _{\infty
}n^{-N\left(  U\right)  } \label{21}
\end{equation}
for some\textbf{\ }$C\left(  \bar{C},V\right)  >0$\textbf{, }while the
functional $N\left(  \cdot\right)  $ is positive for every $U$ smooth.
\end{lemma}

\begin{proof}
\textbf{\ }Our system in the box $\Lambda_{n}$ has $\left(  2n+1\right)  ^{2}$
degrees of freedom, which is hard to study. We are going to fix $\left(
2n+1\right)  ^{2}-\left(  n+1\right)  $ of them, leaving only $n+1$ degrees of
freedom, and we will show that for every choice $\Phi$ of the degrees frozen
we have
\begin{equation}
\left|  \left\langle f\left(  \phi+\psi\right)  |\Phi\right\rangle
_{n,\bar{\phi}}-\left\langle f\left(  \phi\right)  |\Phi\right\rangle
_{n,\bar{\phi}}\right|  \leq C\left(  \bar{C},V\right)  \left|  \left|
f\right|  \right|  _{\infty}n^{-N\left(  U\right)  } \label{22}
\end{equation}
uniformly in $\Phi.$ From that (\ref{21}) evidently follows by integration.
These degrees of freedom are introduced in the following way.

For every $k=0,1,2,...$ we define the layer $L_{k}\subset\mathbb{Z}^{2}$ as
the subset $L_{k}=\left\{  x\in\mathbb{Z}^{2}:\left|  \left|  x\right|
\right|  _{\infty}=k\right\}  .$ For a configuration $\phi$ in $\Lambda_{n}$
we denote by $\Phi_{k},k=0,1,2,...,n$ its restrictions to the layers $L_{k}:$
\[
\Phi_{k}=\phi|_{L_{k}}.
\]
We define now the action $\left(  \psi_{0},\psi_{1},...,\psi_{n}\right)  \phi$
of the group $\left(  \mathbb{S}^{1}\right)  ^{n+1}$ on configurations $\phi$
in $\Lambda_{n}$ by
\[
\left(  \left(  \psi_{0},\psi_{1},...,\psi_{n}\right)  \phi\right)  \left(
x\right)  =\phi\left(  x\right)  +\psi_{k\left(  x\right)  }
\]
where $k\left(  x\right)  =\left|  \left|  x\right|  \right|  _{\infty}$ is
the number of the layer to which the site $x$ belongs. We define the torus
$\Phi\left(  \phi\right)  $ to be the orbit of the configuration $\phi$ under
this action. In other words, $\Phi\left(  \phi\right)  $ is the set of
configurations $\Phi_{0}+\psi_{0},\Phi_{1}+\psi_{1},...,\Phi_{n}+\psi_{n},$
for all possible values of the angles $\psi_{i},$ where the configuration
$\Phi_{k}+\psi_{k}$ on the layer $L_{k}$ is defined by $\left(  \Phi_{k}
+\psi_{k}\right)  \left(  x\right)  =\phi\left(  x\right)  +\psi_{k}.$

Let us fix for every orbit $\Phi$ one representative, $\phi,$ so $\Phi
=\Phi\left(  \phi\right)  ,$ and let $\Phi_{0},\Phi_{1},...,\Phi_{n}$ be the
restrictions, $\Phi_{k}=\phi|_{L_{k}}.$

We will study the conditional Gibbs distribution $\left\langle \cdot
|\Phi\left(  \phi\right)  =\Phi\right\rangle _{n,\bar{\phi}}$. This
distribution is again a Gibbs measure on $\left(  \mathbb{S}^{1}\right)
^{n+1}=\left\{  \left(  \psi_{0},\psi_{1},...,\psi_{n}\right)  \right\}  ,$
corresponding to the nearest neighbour interaction $\mathcal{W}_{\Phi
,\bar{\phi}}=\left\{  W_{k},k=1,2,...,n\right\}  $. It is defined for $k<n$
by
\begin{equation}
W_{k}\left(  \psi_{k},\psi_{k+1}\right)  =\sum_{\substack{x\in L_{k},y\in
L_{k+1}: \\\left|  x-y\right|  =1}}U\left[  \left(  \Phi_{k}+\psi_{k}\right)
\left(  x\right)  ,\left(  \Phi_{k+1}+\psi_{k+1}\right)  \left(  y\right)
\right]  , \label{26}
\end{equation}
while
\begin{equation}
W_{n}\left(  \psi_{n}\right)  =\sum_{\substack{x\in L_{n},y\in L_{n+1}:
\\\left|  x-y\right|  =1}}U\left[  \left(  \Phi_{n}+\psi_{n}\right)  \left(
x\right)  ,\bar{\phi}\left(  y\right)  \right]  . \label{27}
\end{equation}
(Note for the future, that the interactions along the bonds which are
contained within one layer do not contribute to $W$-s.) We are going to show
that for every $k$ the distribution of the random variable $\psi_{k}$ under
$\left\langle \cdot|\Phi\left(  \phi\right)  =\Phi\right\rangle _{n,\bar{\phi
}}$ has a density $p_{k}\left(  t\right)  $ with respect to the Lebesgue
measure on $\mathbb{S}^{1},$ which satisfies
\begin{equation}
\sup_{t\in\mathbb{S}^{1}}\left|  p_{k}\left(  t\right)  -1\right|  \leq
C\sqrt{k}\left(  \frac{n}{k}\right)  ^{-N\left(  U\right)  }, \label{23}
\end{equation}
uniformly in $\Phi,\bar{\phi},$ with $C=C\left(  \bar{C}\right)  .$ That
implies (\ref{22}).

To show (\ref{23}) we note that due to $\mathbb{S}^{1}$-invariance of $U$ we
have
\[
W_{k}\left(  \psi_{k},\psi_{k+1}\right)  =W_{k}\left(  \psi_{k}+\alpha
,\psi_{k+1}+\alpha\right)
\]
for every $\alpha\in\mathbb{S}^{1}.$ Hence $W_{k}\left(  \psi_{k},\psi
_{k+1}\right)  =W_{k}\left(  \psi_{k}-\psi_{k+1},0\right)  ,$ and therefore
the random variables
\[
\chi_{k}=\left\{
\begin{array}
[c]{cc}
\psi_{k}-\psi_{k+1} & \text{ for }k<n\\
\psi_{n} & \text{ for }k=n
\end{array}
\right.
\]
are independent. Since evidently
\begin{equation}
\psi_{k}=\chi_{k}+\chi_{k+1}+...+\chi_{n}, \label{31}
\end{equation}
we are left with the question about the distribution of the sum of independent
random elements of $\mathbb{S}^{1}.$ Were the independent random elements
$\chi_{i}$ identically distributed, with the distribution having density, the
statement (\ref{23}) would be immediate. However, they are not identically
distributed, so we need to work further.

Introducing $W_{k}\left(  \chi_{k}\right)  =W_{k}\left(  \chi_{k},0\right)  $
for $k<n,$ we have that for all $k\leq n$ the distribution of the random
element $\chi_{k}$ is given by the density
\[
q_{k}\left(  t\right)  =\frac{\exp\left\{  -W_{k}\left(  t\right)  \right\}
}{\int\exp\left\{  -W_{k}\left(  t\right)  \right\}  dt}.
\]
Let $t_{\min}$ be (any) global minimum of the function $W_{k}\left(
\cdot\right)  .$ Then for every $t$ the Taylor expansion implies the estimate
\begin{equation}
W_{k}\left(  t_{\min}\right)  \leq W_{k}\left(  t\right)  \leq W_{k}\left(
t_{\min}\right)  +8\bar{C}\left(  k+1\right)  \left|  t-t_{\min}\right|  ^{2},
\label{35}
\end{equation}
due to (\ref{25}), (\ref{26}), (\ref{27}). (This is the point where both
smoothness and two-dimensionality are crucial.) Hence
\begin{equation}
\max q_{k}\left(  t\right)  \leq C_{1}\sqrt{k+1} \label{32}
\end{equation}
for some $C_{1}=C_{1}\left(  \bar{C}\right)  .$

Because of (\ref{31}), $p_{k}\left(  t\right)  =\left(  q_{k}\ast...\ast
q_{n}\right)  \left(  t\right)  ,$ where $\ast$ stays for convolution.
Therefore it is natural to study the Fourier coefficients
\[
a_{s}\left(  q_{l}\right)  =\frac{1}{2\pi}\int_{0}^{2\pi}q_{l}\left(
t\right)  e^{ist}\,dt,
\]
$s=0,\pm1,\pm2,...\,,$ since
\begin{equation}
a_{s}\left(  p_{k}\right)  =\prod_{l=k}^{n}a_{s}\left(  q_{l}\right)  .
\label{61}
\end{equation}
We want to show that for every $s\neq0$\ the last product goes to $0$\ as
$n\rightarrow\infty,$\ uniformly in $s.$ To estimate the coefficients $\left|
a_{s}\left(  q_{l}\right)  \right|  $\ we use the following straightforward

\begin{lemma}
Let $P_{C}$\ be the set of all probability densities $q\left(  \cdot\right)
$\ on a circle, satisfying
\[
\sup_{t\in\mathbb{S}^{1}}q\left(  t\right)  \leq C,
\]
and $s$\ be an integer. Then the functional on $P_{C},$\ given by the
integral
\[
\frac{1}{2\pi}\int_{0}^{2\pi}q\left(  t\right)  \cos\left(  st\right)  \,dt,
\]
attains its maximal value at the density
\[
q_{C}\left(  t\right)  =\left\{
\begin{array}
[c]{ll}
C & \text{ if }\left|  t-\frac{2\pi k}{s}\right|  \leq\frac{1}{2Cs}\text{ for
some }k=0,...,s-1,\\
0 & \text{otherwice.}
\end{array}
\right.
\]
\end{lemma}

Using this lemma and the estimate (\ref{32}), we obtain that
\begin{align}
&  \sup\left\{  \left|  a_{s}\left(  q_{l}\right)  \right|  :s\neq0\right\}
\label{33}\\
&  \leq2C_{1}\sqrt{l+1}\int_{0}^{\frac{1}{2}\left(  C_{1}\sqrt{l+1}\right)
^{-1}}\left(  1-\frac{t^{2}}{3}\right)  dt=1-\frac{1}{36\left(  C_{1}\right)
^{2}\left(  l+1\right)  }.\nonumber
\end{align}
Since
\[
\sup_{t\in\mathbb{S}^{1}}\left|  p_{k}\left(  t\right)  -1\right|  \leq
\sum_{s\neq0}\left|  a_{s}\left(  p_{k}\right)  \right|  ,
\]
we are almost done. Namely, note that due to Parseval identity and (\ref{32})
we have for every $l$\
\[
1+\sum_{s\neq0}\left|  a_{s}\left(  q_{l}\right)  \right|  ^{2}=\int\left(
q_{l}\left(  t\right)  \right)  ^{2}dt\leq C_{1}\sqrt{l+1}.
\]
Let us introduce now the densities\textbf{\ }$p_{k,r}\left(  t\right)
=\left(  q_{k}\ast...\ast q_{r}\right)  \left(  t\right)  ,\,k\leq r\leq
n.$\textbf{\ }Due to Cauchy inequality,
\[
1+\sum_{s\neq0}\left|  a_{s}\left(  p_{k,k+1}\right)  \right|  \leq
C_{1}\sqrt[4]{\left(  k+1\right)  \left(  k+2\right)  }.
\]
Therefore by (\ref{33}) and (\ref{61})
\begin{equation}
\sup_{t\in\mathbb{S}^{1}}\left|  p_{k,r}\left(  t\right)  -1\right|  \leq
C_{1}\sqrt[4]{\left(  k+1\right)  \left(  k+2\right)  }\prod_{l=k+2}
^{r}\left(  1-\frac{1}{36\left(  C_{1}\right)  ^{2}\left(  l+1\right)
}\right)  , \label{83}
\end{equation}
which ends the proof of (\ref{23}), with $C=2C_{1}\left(  \bar{C}\right)
$\ and $N\left(  U\right)  =\frac{1}{36\left(  C_{1}\left(  \bar{C}\right)
\right)  ^{2}}.$
\end{proof}

\subsection{Theorem 1: Singular case.}

The key step in the above proof was the use of the Taylor expansion, to bound
the densities $q_{r}.$\textbf{\ }There the existence of the second derivative
of $U$\ and its boundedness was used in a crucial way. Yet, one can use
essentially the same arguments to treat the general case, without smoothness
assumption. The main idea is to represent the singular interaction as a small
perturbation of a smooth one, smallness being understood in the $L_{1}$ sense.
Another version of this idea was used earlier in \cite{BI,BCPK,DV,IV}.

Namely, we will consider the nearest neighbour interaction
\begin{equation}
\bar{U}\left(  \phi\right)  =U\left(  \phi\right)  -\upsilon\left(
\phi\right)  , \label{41}
\end{equation}
where $U$ is a smooth function with a bounded second derivative, as above,
while $\upsilon\geq0$ is a ''small'' singular component. The precise meaning
of smallness will be made explicit a bit later, see (\ref{51}). However,
already now we can say that every continuous function $\bar{U}$ can be written
in the form (\ref{41}), with $U$ twice differentiable and with $\upsilon$
satisfying
\begin{equation}
0\leq\upsilon\left(  \cdot\right)  \leq\varepsilon, \label{42}
\end{equation}
with $\varepsilon>0$ arbitrarily small. That follows immediately for example
from the Weierstrass theorem, stating that the trigonometric polynomials are
everywhere dense in the space of continuous functions on the circle. Clearly,
the estimate (\ref{42}) implies $L_{1}$-smallness of $\upsilon$, whatever the
latter may mean.

We will denote by $\mathcal{\bar{H}}$\textbf{\ }the Hamiltonian corresponding
to the singular interaction $\bar{U},$\ while\textbf{\ }$\mathcal{H}
$\textbf{\ }will be the Hamiltonian defined by the smooth part of interaction,
$U.$\ To proceed with the expansion, we introduce the set\textbf{\ }
$\mathcal{E}_{n}$\textbf{\ }to be the collection of all bonds of\textbf{\ }
$\mathbb{Z}^{2}$ \ with at least one end in the box $\Lambda_{n},$\ and
rewrite the partition function $Z_{n}^{\bar{U},\bar{\phi}}$ in $\Lambda_{n},$
corresponding to the interaction $\bar{U}$ and the boundary conditions
$\bar{\phi}$, as follows:
\begin{align*}
Z_{n}^{\bar{U},\bar{\phi}}  &  =\int_{\Omega_{n}}\exp\left\{  -\mathcal{\bar
{H}}\left(  \mathbf{\phi|}\bar{\phi}\right)  \right\}  d\mathbf{\phi}\\
\  &  =\int_{\Omega_{n}}\exp\left\{  -\mathcal{H}\left(  \mathbf{\phi|}
\bar{\phi}\right)  \right\}  \prod_{\left\langle x,y\right\rangle
\in\mathcal{E}_{n}}\left[  1+\left(  e^{\upsilon\left(  \mathbf{\phi}\left(
x\right)  -\mathbf{\phi}\left(  y\right)  \right)  }-1\right)  \right]
d\mathbf{\phi}\\
\  &  =\sum_{A\subset\mathcal{E}_{n}}\int_{\Omega_{n}}\exp\left\{
-\mathcal{H}\left(  \mathbf{\phi|}\bar{\phi}\right)  \right\}  \prod
_{\left\langle x,y\right\rangle \in A}\left(  e^{\upsilon\left(  \mathbf{\phi
}\left(  x\right)  -\mathbf{\phi}\left(  y\right)  \right)  }-1\right)
d\mathbf{\phi}\\
\  &  \equiv\sum_{A\subset\mathcal{E}_{n}}Z_{n}^{U,\bar{\phi},A}.
\end{align*}
For every subset $A\subset\mathcal{E}_{n}$ we now introduce the probability
distribution
\[
\mu_{n}^{U,\bar{\phi},A}\left(  d\mathbf{\phi}\right)  =\frac{1}{Z_{n}
^{U,\bar{\phi},A}}\exp\left\{  -\mathcal{H}\left(  \mathbf{\phi|}\bar{\phi
}\right)  \right\}  \prod_{\left\langle x,y\right\rangle \in A}\left(
e^{\upsilon\left(  \mathbf{\phi}\left(  x\right)  -\mathbf{\phi}\left(
y\right)  \right)  }-1\right)  d\mathbf{\phi.}
\]
Then we have for the original Gibbs state $\mu_{n}^{\bar{U},\bar{\phi}}$ the
following decomposition:
\[
\mu_{n}^{\bar{U},\bar{\phi}}=\sum_{A\subset\mathcal{E}_{n}}\pi_{n}\left(
A\right)  \mu_{n}^{U,\bar{\phi},A},
\]
with the probabilities $\pi_{n}\left(  \cdot\right)  $ given by
\[
\pi_{n}\left(  A\right)  =\frac{Z_{n}^{U,\bar{\phi},A}}{Z_{n}^{\bar{U}
,\bar{\phi}}}.
\]
Note that the states $\mu_{n}^{U,\bar{\phi},A}$ are themselves Gibbs states in
$\Lambda_{n}$, corresponding to the boundary condition $\bar{\phi}$ and the
(non-translation invariant) nearest neighbour interaction $\mathcal{U}^{A},$
which for bonds outside $A$ is given by our smooth function $U\left(  \phi
_{s}-\phi_{t}\right)  ,$ while on bonds from $A$ it equals to $U\left(
\phi_{s}-\phi_{t}\right)  -\ln\left(  e^{\upsilon\left[  \phi_{s}-\phi
_{t}\right]  }-1\right)  .$ (Here the positivity of the function $\upsilon$ is
used.) Let us now introduce the bond percolation process $\mathcal{A}$ on
$\mathcal{E}_{n},$ defining its probability distribution $\mathbb{P}_{n}$ by
\[
\mathbb{P}_{n}\left(  \mathcal{A}=A\right)  =\pi_{n}\left(  A\right)  .
\]
This process is of course a dependent percolation process. Happily, it turns
out that it is dominated by independent bond percolation, with probability of
a bond to be open very small! Our claim would follow once we check that the
conditional probabilities
\[
\mathbb{P}_{n}\left(  b\in\mathcal{A}|\left(  \mathcal{E}_{n}\backslash
b\right)  \cap\mathcal{A}=\mathcal{D}\right)
\]
are small uniformly in $\mathcal{D}.$ We will show this under the following
condition on the smallness of the singular part $\upsilon$ of the interaction
$\bar{U}.$ We suppose that

\begin{itemize}
\item $\bar U\left(  \phi\right)  =U\left(  \phi\right)  -\upsilon\left(
\phi\right)  ,$ with $U$ having bounded second derivative,

\item $\upsilon\ge0,$

\item  for every choice of the four values $\phi_{1},\phi_{2},\phi_{3}
,\phi_{4}$
\begin{equation}
\frac{\int\exp\left\{  -\sum_{i=1}^{4}U\left(  \phi-\phi_{i}\right)
+\sum_{i=1}^{4}\upsilon\left(  \phi-\phi_{i}\right)  \right\}  d\phi}{\int
\exp\left\{  -\sum_{i=1}^{4}U\left(  \phi-\phi_{i}\right)  \right\}  d\phi
}\leq1+\varepsilon, \label{51}
\end{equation}
with $\varepsilon$ small enough.
\end{itemize}

In words, the last condition says that the expectation of the observable
$\exp\left\{  \sum_{i=1}^{4}\upsilon\left(  \phi-\phi_{i}\right)  \right\}  $
with respect to a single site conditional Gibbs distribution corresponding to
the (smooth) interaction $U$ and any boundary condition $\phi_{1},\phi
_{2},\phi_{3},\phi_{4}$ around that site, is smaller than $1+\varepsilon.$ A
straightforward calculation implies that under (\ref{51})
\begin{equation}
\mathbb{P}_{n}\left(  b\in\mathcal{A}|\left(  \mathcal{E}_{n}\backslash
b\right)  \cap\mathcal{A}=\mathcal{D}\right)  \leq\varepsilon, \label{LSSS}
\end{equation}
uniformly in $\mathcal{D}.$ We denote by $\mathbb{Q}_{\varepsilon}$ the
distribution of the corresponding independent bond percolation process,
$\eta_{\cdot}$.

The strategy of the remainder of this subsection is the following:

\begin{itemize}
\item  we will show that if the set $A$ is sparse enough, then for the measure
$\mu_{n}^{U,\bar{\phi},A}$ the analog of the estimate (\ref{21}) holds.

\item  such sparse sets $A$ constitute the dominant contribution to the
distribution $\mathbb{P}_{n}.$
\end{itemize}

Let us formulate now the sparseness condition on $A$ we need.

In what follows, by a path we will mean a sequence of pairwise distinct bonds
of our lattice, such that any two consecutive bonds share a site. A path with
coinciding beginning and end is called a loop. If a loop surrounds the origin,
we will call it a circuit. Any two objects of the above will be called
disjoint, if they share neither a bond nor a site. The same objects,
associated with the dual lattice will be called d-sites, d-bonds, d-paths,
d-loops and d-circuits.

Suppose the set $A$ is given, and $\lambda_{1},\lambda_{2,}...,\lambda_{\nu}$
be a collection of disjoint d-circuits, avoiding $A.$ The latter means that no
d-bond of any $\lambda_{k}$ crosses any of the bonds from $A.$ We suppose that
these d-circuits are ordered by ``inclusion''. Then we introduce layers
$L_{k}$ by
\[
L_{k}=\left\{  x\in\mathbb{Z}^{2}:x\in\mathrm{Int}\left(  \lambda_{k}\right)
\backslash\mathrm{Int}\left(  \lambda_{k-1}\right)  \right\}  ,\,k=1,2,...,\nu
+1,
\]
with the convention that $\mathrm{Int}\left(  \lambda_{0}\right)  =\emptyset$
and $\mathrm{Int}\left(  \lambda_{\nu+1}\right)  =\mathbb{Z}^{2}.$ (Note that
these layers are connected sets of sites, and they surround the origin in the
same way as the ''old'' layers did.) For every configuration $\phi$ in
$\Lambda_{n}$ we introduce, as in the previous section, the layer
configurations $\Phi_{k},k=1,2,...,\nu+1$ as its restrictions to the layers
$L_{k},$ the layer angles $\psi_{1},...,\psi_{\nu},$ the $\nu$-dimensional
torus $\Phi\left(  \phi\right)  ,$ and we note that the distribution of $\psi
$-s under the condition that the orbit $\Phi\left(  \phi\right)  $ is fixed,
is a (one-dimensional) Gibbs distribution. Moreover, it is defined by the
nearest neighbour interaction $\mathcal{W}_{\Phi,\bar{\phi}}=\left\{
W_{k},k=1,2,...,\nu\right\}  $, given by almost the same formula, as
(\ref{26}): for $k<\nu$
\begin{equation}
W_{k}\left(  \psi_{k},\psi_{k+1}\right)  =\sum_{\substack{x\in L_{k},y\in
L_{k+1}: \\\left|  x-y\right|  =1}}U\left[  \left(  \Phi_{k}+\psi_{k}\right)
\left(  x\right)  ,\left(  \Phi_{k+1}+\psi_{k+1}\right)  \left(  y\right)
\right]  , \label{56}
\end{equation}
while for $k=\nu$
\begin{equation}
W_{\nu}\left(  \psi_{\nu}\right)  =\sum_{\substack{x\in L_{\nu},y\in L_{\nu
+1}: \\\left|  x-y\right|  =1}}U\left[  \left(  \Phi_{n}+\psi_{n}\right)
\left(  x\right)  ,\left(  \phi\vee\bar{\phi}\right)  \left(  y\right)
\right]  . \label{57}
\end{equation}
(Here the configuration $\phi\vee\bar{\phi}$ equals to $\phi$ inside
$\Lambda_{n}$ and to $\bar{\phi}$ outside $\Lambda_{n}.$) Note that the
singular part of the interaction $\mathcal{U}^{A}$ does not enter in these
formulas, precisely because the d-circuits $\lambda_{k}$ avoid the set $A$!
Hence we can conclude that for every $k$ the distribution of the random
variable $\psi_{k}$ under the measure $\left\langle \cdot|\Phi\left(
\phi\right)  =\Phi\right\rangle _{n,\bar{\phi}}$ has a density $p_{k}\left(
t\right)  $ on $\mathbb{S}^{1},$ which satisfies the following analog of
(\ref{83}):
\begin{equation}
\sup_{t\in\mathbb{S}^{1}}\left|  p_{k}\left(  t\right)  -1\right|  \leq
C_{1}\sqrt[4]{\left|  \lambda_{k}\right|  \left|  \lambda_{k+1}\right|  }
\exp\left\{  -\frac{1}{36\left(  C_{1}\right)  ^{2}}\sum_{l=k+2}^{\nu}\frac
{1}{\left|  \lambda_{l}\right|  }\right\}  , \label{53}
\end{equation}
uniformly in $\Phi,\bar{\phi}.$ The last relation suggests the following

\begin{definition}
\textbf{of sparseness:} The set $A$ of bonds in $\mathcal{E}_{n}$ is $\tau
$-sparse, if there exists a family of $\nu\left(  A\right)  $ disjoint
d-circuits $\lambda_{l}$ in $\Lambda_{n}$, avoiding $A,$ and such that
\[
\sum_{l=1}^{\nu\left(  A\right)  }\frac{1}{\left|  \lambda_{l}\right|  }
\geq\tau\ln n.
\]
\end{definition}

Therefore we will be done, once we show the following:

\begin{proposition}
\label{tau}For any $\kappa,$ $1>\kappa>0,$ there exists a value $\tau
=\tau\left(  \kappa\right)  >0,$ such that
\begin{equation}
\mathbb{P}_{n}\left(  \mathcal{A}\text{ is not }\tau\text{-sparse}\right)
\leq e^{-n^{\kappa}}. \label{60}
\end{equation}
\end{proposition}

The proof of this proposition is the content of the following subsections.

\subsubsection{$\tau$-sparseness is typical.}

For every $l=2,3,...$ let us define the northern rectangle
\[
R_{N}^{l}=[-2^{l},\dots,2^{l}]\times\lbrack2^{l-1}+1,\dots,2^{l}],
\]
and let the eastern, southern and western rectangles $R_{E}^{l},\,R_{S}^{l}$
and $R_{W}^{l}$ be the clock-wise rotations of $R_{N}^{l}$ by, respectively,
$\pi/2$, $\pi$ and $3\pi/2$ with respect to the origin. Define the $l$-th
shell $T^{l}$ by
\[
T^{l}~=~R_{N}^{l}\cup R_{E}^{l}\cup R_{S}^{l}\cup R_{W}^{l}.
\]
Clearly, $T^{l}\subset\Lambda_{n}$ once $n\geq2^{l},$ while different $T^{l}
$-s are disjoint.

Let a configuration $A$ of bonds be given. By a good crossing of a rectangle
$R_{\cdot}^{\cdot}$ we will mean a d-path, joining the two short sides of
$R_{\cdot}^{\cdot}$ and avoiding $A.$ We denote the set of such crossings by
$\mathcal{R}^{\longleftrightarrow}.$ Let $\lambda_{N}^{l},\lambda_{E}
^{l},\lambda_{S}^{l},\lambda_{W}^{l}$ be four good crossings of the rectangles
$R_{N}^{l},R_{E}^{l},R_{S}^{l},R_{W}^{l}$ respectively. Then the collection of
those d-bonds of the union $\lambda_{N}^{l}\cup\lambda_{E}^{l}\cup\lambda
_{S}^{l}\cup\lambda_{W}^{l},$ which are seen from the origin, form a d-circuit
avoiding $A.$ Therefore we want to get a

\subsubsection{Lower bound on the number of disjoint good crossing of a rectangle}

We claim that for all $\varepsilon$ sufficiently small there exist
$\alpha=\alpha(\varepsilon)>0$ and $c_{1}=c_{1}(\varepsilon)>0$ such that at
each scale $k$ the $\mathbb{Q}_{\varepsilon}$-probability that there are less
than $\mathbb{\alpha}2^{k}$ disjoint good crossings of $R_{N}^{k}$ is smaller
than $e^{-c_{1}2^{k}}$, where $\mathbb{Q}_{\varepsilon}$ is the measure of the
independent bond percolation process $\eta_{\cdot},$ defined after~(\ref{LSSS}).

Indeed, by the Ford-Fulkerson min-cut/max-flow Theorem (see e.g. \cite{R}),
the number of disjoint good crossings of $R_{N}^{k}$ (which by definition are
left-to-right crossings by d-paths) is bounded from below by
\[
\frac{1}{2}\min_{{\widetilde{\lambda}\in}\mathcal{R}^{\updownarrow}}\left\{
\left|  {\widetilde{\lambda}}\right|  -\left|  {\widetilde{\lambda}}\cap
A\right|  \right\}  ,
\]
where the minimum is taken over the set $\mathcal{R}^{\updownarrow}$ of all
''cuts'', which are just paths in {$R_{N}^{k},$ joining the bottom and top
sides of $R_{N}^{k}.$ The min-cut quantity }$\min_{{\widetilde{\lambda}}
}\left\{  \left|  {\widetilde{\lambda}}\right|  -\left|  {\widetilde{\lambda}
}\cap A\right|  \right\}  $ equals to the maximal left-to-right flow by
d-paths, avoiding $A,$ and the factor $1/2$ accounts for the fact that the
corresponding d-paths might share the same d-sites, so in order to estimate
the number of disjoint paths we have to take a half of the total flow.

Evidently,
\begin{equation}
\mathbb{Q}_{\varepsilon}\left(  \exists{\widetilde{\lambda}\in}\mathcal{R}
^{\updownarrow}\text{ with }\left|  {\widetilde{\lambda}}\right|  -\left|
{\widetilde{\lambda}}\cap A\right|  \leq\alpha\text{$2^{k}$}\right)  \leq
\sum_{{\widetilde{\lambda}\in}\mathcal{R}^{\updownarrow}}\mathbb{Q}
_{\varepsilon}\left(  \left|  {\widetilde{\lambda}}\right|  -\left|
{\widetilde{\lambda}}\cap A\right|  \leq\alpha\text{$2^{k}$}\right)  ,
\label{gr}
\end{equation}
while for every ${\widetilde{\lambda}}$
\[
\mathbb{Q}_{\varepsilon}\left(  \left|  {\widetilde{\lambda}}\right|  -\left|
{\widetilde{\lambda}}\cap A\right|  \leq\alpha\text{$2^{k}$}\right)
\leq2^{\left|  {\widetilde{\lambda}}\right|  }\varepsilon^{\left|
{\widetilde{\lambda}}\right|  -\mathbb{\alpha}2^{k}}~\leq~e^{-c_{2}\left|
{\widetilde{\lambda}}\right|  },
\]
since any top-to-bottom crossing contains at least $2^{k-1}$ bonds. Here
$c_{2}=c_{2}(\alpha,\varepsilon)>0$ satisfies
\[
\lim_{\varepsilon\rightarrow0}c_{2}(\alpha,\varepsilon)~=~\infty,
\]
once $\alpha<1/2$. Thus, choosing $\alpha<1/2$ and $\varepsilon$ sufficiently
small, we infer that there exists $c_{1}>0$, such that the right hand side of
(\ref{gr}) is bounded above by
\[
2^{k}\sum_{l=2^{k-1}}^{\infty}3^{l}e^{-c_{2}l}~\leq~~e^{-c_{1}2^{k}}.
\]
Thus, the min-cut/max-flow theorem insures that up to the $\mathbb{Q}
_{\varepsilon}$-probability $1-e^{-c_{1}2^{k}}$, there are at least
$\alpha2^{k-1}$ disjoint good crossings $\lambda_{i}$ of $R_{N}^{k}$.
Moreover, observe that at least $\alpha2^{k-2}$ of these d-paths have the
length bounded above by $\alpha^{-1}2^{k+3}$. Indeed, should this not be the
case,
\[
\sum_{i}|\lambda_{i}|>\alpha2^{k-2}\,\frac{1}{\alpha}2^{k+3}=2|R_{N}^{k}|
\]
which in view of the disjointedness of $\lambda_{i}$-s is impossible.

Let us say that a left-to-right crossing d-path $\lambda$ of the $k$-th scale
is $\alpha$-short, if $|\lambda|<\alpha^{-1}2^{k+3},$ and define the event
\[
\mathcal{T}_{N}^{k,\alpha}~=~\left\{  A:{\text{there are at least }
\alpha2^{k-2}\text{ disjoint good $\alpha$-short crossings of $R_{N}^{k}$}
}\right\}  {\text{$.$ }}
\]
What we have proved up to now can be summarized as follows:

There exists $c_{1}>0$, such that uniformly in $k$,
\begin{equation}
\mathbb{Q}_{\varepsilon}\left(  \mathcal{T}_{N}^{k,\alpha}\right)
~\geq1-~e^{-c_{1}2^{k}}, \label{short_crossings}
\end{equation}
as soon as $\alpha$ and $\varepsilon$ are sufficiently small.

\subsubsection{Proof of Proposition~\ref{tau}}

Consider now the event
\[
\mathcal{T}^{k,\alpha}=\mathcal{T}_{N}^{k,\alpha}\cap\mathcal{T}_{E}
^{k,\alpha}\cap\mathcal{T}_{S}^{k,\alpha}\cap\mathcal{T}_{W}^{k,\alpha}.
\]
From the previous argument one knows that for $\varepsilon$ close enough to
$0$ the $\mathbb{Q}_{\varepsilon}$-probability of the event $\mathcal{T}
^{k,\alpha}$ is at least $1-4\,e^{-c_{1}2^{k}}$. Note that under
$\mathcal{T}^{k,\alpha}$ there are at least $\alpha2^{k-2}$ disjoint
d-circuits in $T^{k}$, avoiding $A,$ all of which have length at most
$2^{k+5}/\alpha$. Also, the events $\mathcal{T}^{k,\alpha}$ are
non-decreasing, therefore their $\mathbb{P}_{n}$-probability is at least
$1-4\,e^{-c_{1}2^{k}}$ as well.

The claim of Proposition~\ref{tau} is now an immediate consequence: Let
$1>\rho>0$. Then, for every $n=2,3,\dots$ the event
\[
\mathcal{T}^{n,\rho,\alpha}=\bigcap_{k=\left[  \rho\log_{2}\,n\right]
}^{\left[  \log_{2}\,n\right]  }\mathcal{T}^{k,\alpha}
\]
has, by (\ref{short_crossings}), $\mathbb{P}_{n}$-probability at least
$1-c_{3}e^{-c_{4}n^{\rho}}$. However, by the very construction, the occurrence
of the event $\mathcal{T}^{n,\rho,\alpha}$ ensures that in each shell $T_{k}$,
$k\in\{\left[  \rho\log_{2}\,n\right]  ,\dots,\left[  \log_{2}\,n\right]  \}$,
it is possible to find a family of disjoint d-circuits avoiding $A$ and such
that the sum of the inverse of their lengths is at least $\alpha^{2}/128$.
Their total is at least
\[
\frac{\alpha^{2}}{128}\frac{1-\rho}{2}\log_{2}\,n.
\]
The conclusion (\ref{60}) follows.

\subsubsection{General finite-range interactions}

We briefly describe the main modifications to the proof given above, which are
needed in order to treat the case of finite-range, non nearest-neighbour
interactions $\bar{U}_{\Lambda}$, $\Lambda\Subset{\mathbb{Z}}^{2}$. As
in~\eqref{41}, we decompose $\bar{U}_{\Lambda}= U_{\Lambda}- \upsilon
_{\Lambda}$ to a smooth part $U_{\Lambda}$ and a small singular part
$0\leq\upsilon_{\Lambda}\leq\varepsilon$. Notice that the choice of
$\varepsilon= \varepsilon(r_{\Lambda})$ will in general depend on the diameter
$r_{\Lambda}=\mathrm{diam} (\Lambda)$ of the interaction set $\Lambda$.

The singular part of the interaction will be controlled by a dependent
\emph{site} percolation process, which we construct in two steps as follows.
Define $\bar\Lambda_{n} = \{ x:\, x+\Lambda\cap\Lambda_{n} \neq\emptyset\}$.

\noindent\textbf{Step 1.} As in the nearest-neighbour case, write
\begin{align*}
Z_{n}^{\bar{U},\bar{\phi}}  &  =\sum_{A\subset\bar{\Lambda}_{n}}\int
_{\Omega_{n}}\exp\left\{  -\mathcal{H}\left(  \mathbf{\phi|}\bar{\phi}\right)
\right\}  \prod_{x\in A}\left(  e^{\upsilon_{\Lambda}\left(  \mathbf{\phi
}_{\cdot\,+x}\right)  }-1\right)  d\mathbf{\phi}\\
&  \overset{\triangle}{=}\sum_{A\subset\bar{\Lambda}_{n}}Z_{n}^{U,\bar{\phi
},A}\,.
\end{align*}
Then, exactly as before, it is easy to show that the probability distribution
\begin{equation}
\mathbb{P}_{n}\left(  \mathcal{A}=A\right)  \overset{\triangle}{=}\frac
{Z_{n}^{U,\bar{\phi},A}}{Z_{n}^{\bar{U},\bar{\phi}}} \label{PdeA}
\end{equation}
on $\{0,1\}^{\bar{\Lambda}_{n}}$ is stochastically dominated by the Bernoulli
site percolation process ${\mathbb{Q}}_{\varepsilon}$ with density
$\varepsilon$.

\noindent\textbf{Step 2.} Let us split ${\mathbb{Z}}^{2}$ into the disjoint
union of the shifts of squares $B_{\Lambda}\overset{\triangle}{=}
\{-2r_{\Lambda},\dots,2r_{\Lambda}\}^{2}$,
\[
{\mathbb{Z}}^{2}=\bigvee_{x}\left(  4r_{\Lambda}\,x+B_{\Lambda}\right)  \,.
\]
Given a realization $A$ of the random set ${\mathcal{A}}$ (distributed
according to~(\ref{PdeA})) let us say that $x\in{\mathbb{Z}}^{2}$ is good if
$4r_{\Lambda}\,x+B_{\Lambda}\cap A=\emptyset$. Thus, for every $n$,
${\mathcal{A}}$ induces a probability distribution on $\{0,1\}^{{\mathbb{Z}
}^{2}}$, which stochastically dominates Bernoulli site percolation with
density $1-(1-\varepsilon)^{16r_{\Lambda}^{2}}$.

This dictates the choice of $\varepsilon$ in terms of the diameter of the
interaction $r_{\Lambda}$: For example, $\varepsilon= 1/(Cr_{\Lambda}^{2})$
for $C$ large enough qualifies.

The end of the proof is a straightforward modification of the one in the
nearest-neighbour case.

\subsection{Long-range case: Proof of Theorem~\ref{thm2}}

In this section we study the long-range case, by adapting the technique of~
\cite{P, FP} to the setting of singular interaction. As in the previous
section, we restrict our attention to the case of $\mathbb{S}^{1}$-valued
spins (the extension to the general case is done in the same way as before).
We give here a proof only for the case when all the interactions $J_{x}$
in\textbf{ }(\ref{70}) are nonnegative. The proof in the general case is then straightforward.

Let again $\Lambda_{n}$ be the box $\left\{  x\in\mathbb{Z}^{2}:\left|
\left|  x\right|  \right|  _{\infty}\leq n\right\}  $, $\mathcal{E}_{n}
^{J}=\{\{x,y\}\,:\,J_{x-y}\neq0,\{x,y\}\cap\Lambda_{n}\neq\emptyset\}$, and
let $\bar{\phi}$ be an arbitrary boundary condition outside $\Lambda_{n}$. The
relative Hamiltonian takes the form
\[
\mathcal{\bar{H}}\left(  \mathbf{\phi}_{\Lambda_{n}}|\bar{\phi}\right)
=\sum_{\substack{\{x,y\}\in\mathcal{E}_{n}^{J} \\\{x,y\}\subset\Lambda_{n}
}}J_{x-y}\,\bar{U}(\phi_{x}-\phi_{y})+\sum_{\substack{\{x,y\}\in
\mathcal{E}_{n}^{J} \\\{x,y\}\not \subset\Lambda_{n}}}J_{x-y}\,\bar{U}
(\phi_{x}-\bar{\phi}_{y})\,,
\]
where as in (\ref{41}) the interaction $\bar{U}$ consists of smooth part $U$
and small part $\upsilon$.

Recall that due to the normalization assumption (\ref{J_normalized}),\textbf{
} we can interpret the numbers $j(x)\overset{\triangle}{=}J_{x}$ as the
transition probabilities of a symmetric random-walk $X_{\cdot}$ on
${\mathbb{Z}}^{2}$. We denote by ${\mathbb{E}}_{X}$ expectation w.r.t. this
random-walk conditioned to start at the origin at time $0$. Our assumption on
the coupling constants $J_{\cdot}$\ is that $X_{\cdot}$\ is recurrent.

Let $\left\langle \cdot\right\rangle _{n,\bar\phi}$ be the Gibbs state in
$\Lambda_{n}$ corresponding to the interaction $\bar U$ and the boundary
condition $\bar\phi.$ To prove the theorem, it is enough to show that, for any
$\delta>0$, any bounded local function $f\left(  \phi\right)  $ and any
$\psi\in\mathbb{S}^{1}$,
\begin{equation}
\lim_{n\to\infty}\left|  \left\langle f\left(  \phi+\psi\right)  \right\rangle
_{n,\bar{\phi} }-\left\langle f\left(  \phi\right)  \right\rangle
_{n,\bar{\phi}}\right|  \leq\delta\,. \label{decaylongrange}
\end{equation}

\subsubsection{Expansion of the measure}

As in Subsection~2.2, we expand the Gibbs measure as
\[
\mu_{n}^{\bar{U},\bar{\phi}}=\sum_{A\subset\mathcal{E}_{n}^{J}}\pi_{n}\left(
A\right)  \mu_{n}^{U,\bar{\phi},A},
\]
with the probabilities $\pi_{n}\left(  \cdot\right)  $ given by
\[
\pi_{n}\left(  A\right)  =\frac{Z_{n}^{U,\bar{\phi},A}}{Z_{n}^{\bar{U}
,\bar{\phi}}},
\]
and consider the bond percolation process $\mathcal{A}$ on $\mathcal{E}
_{n}^{J}$ with probability distribution
\[
\mathbb{P}_{n}\left(  \mathcal{A}=A\right)  =\pi_{n}\left(  A\right)  .
\]
Exactly as before, we can show that this process is stochastically dominated
by independent bond percolation process $\mathbb{Q}_{J,\varepsilon}$ on
$\mathcal{E}_{n}^{J}$ with probabilities
\[
\mathbb{Q}_{J,\varepsilon}(\{x,y\}\in\mathcal{A})=\varepsilon J_{x-y}.
\]
From now on, we always assume that $\varepsilon$ is chosen strictly smaller
than $1$.\newline We will use the following notation for the connectivities of
the process $\mathbb{Q}_{J,\varepsilon}$:
\[
p_{x,\varepsilon}=\mathbb{Q}_{J,\varepsilon}\left(  0\overset{A}
{\leftrightarrow}x\right)  \,.
\]
Notice that
\[
p_{x,\varepsilon}\leq\sum_{n=1}^{\infty}\varepsilon^{n}j^{(n)}(x)\overset
{\triangle}{=}d_{\varepsilon}(x)\,,
\]
where $j^{(n)}$ are the $n$-steps transition probabilities of the random-walk
$X_{\cdot}$. Therefore
\begin{equation}
c(\varepsilon)\overset{\triangle}{=}\sum_{x}p_{x,\varepsilon}\leq\sum
_{x}d_{\varepsilon}(x)=\sum_{n=1}^{\infty}\varepsilon^{n}=\frac{\varepsilon
}{1-\varepsilon}\,, \label{cepsilon}
\end{equation}
and the numbers $c(\varepsilon)^{-1}\,p_{x,\varepsilon}$ can be considered as
the transition probabilities of a new random-walk on ${\mathbb{Z}}^{2}$, which
we denote by $Y_{\cdot}$; expectation w.r.t. $Y_{\cdot}$ conditioned to start
at $0$ at time $0$ is denoted by ${\mathbb{E}}_{Y}$. The following lemma plays
an essential role in the sequel:

\begin{lemma}
$X_{\cdot}$ recurrent $\Longrightarrow\,\,Y_{\cdot}$ recurrent. \label{Yrecurrent}
\end{lemma}

\begin{proof}
The recurrence of $X_{\cdot}$ is equivalent (see Th. 8.2 in Chapter~II of
\cite{Sp}) to
\begin{equation}
\int_{{\mathbb{T}}^{2}}\frac{d\theta}{1-\phi(\theta)}=\infty\,,
\label{Xdiverges}
\end{equation}
where
\[
\phi(\theta)\overset{\triangle}{=}{\mathbb{E}}_{X}e^{i(\theta,X_{1})}=\sum
_{x}e^{i(\theta,x)}\,j(x)=\sum_{x}\cos\left(  (\theta,x)\right)  \,j(x)\,,
\]
One has to show that
\begin{equation}
\int_{{\mathbb{T}}^{2}}\frac{d\theta}{1-{\mathbb{E}}_{Y}e^{i(\theta,Y_{1})}
}=\infty\,. \label{Ydiverges}
\end{equation}
Now, $Y_{\cdot}$ is symmetric. Thus
\begin{align*}
1-{\mathbb{E}}_{Y}e^{i(\theta,Y_{1})}  &  ={\mathbb{E}}_{Y}\left(
1-\cos\left(  (\theta,Y_{1})\right)  \right) \\
&  =\tfrac{1}{c(\varepsilon)}\sum_{x}\left(  1-\cos\left(  (\theta,x)\right)
\right)  \,p_{x,\varepsilon}\\
&  \leq\tfrac{1}{c(\varepsilon)}\sum_{x}\sum_{n=1}^{\infty}\left(
1-\cos\left(  (\theta,x)\right)  \right)  \,\varepsilon^{n}j^{(n)}(x)\\
&  =\tfrac{1}{c(\varepsilon)}\sum_{n=1}^{\infty}\left(  1-\phi^{n}
(\theta)\right)  \,\varepsilon^{n}\\
&  =\frac{1-\phi(\theta)}{c(\varepsilon)}\,\sum_{n=1}^{\infty}\varepsilon
^{n}\,\left(  1+\phi(\theta)+\dots+\phi^{n-1}(\theta)\right) \\
&  =\frac{1-\phi(\theta)}{c(\varepsilon)}\,\sum_{n=0}^{\infty}\phi^{n}
(\theta)\,\sum_{k>n}\varepsilon^{k}\\
&  =\frac{(1-\phi(\theta))\varepsilon}{c(\varepsilon)(1-\varepsilon
)(1-\varepsilon\phi(\theta))}\,,
\end{align*}
which implies that~(\ref{Ydiverges}) follows from~(\ref{Xdiverges}).
\end{proof}

\subsubsection{The spin-wave}

Let us denote by $V$ the support of $f$ .

Given a subset $A\subseteq\mathcal{E}_{n}^{J},$ we define the equivalence
relation $\overset{A}{\leftrightarrow}$ between sites of $\mathbb{Z}^{2}$ by
saying that $x\overset{A}{\leftrightarrow}y$ iff there is a path made from the
bonds of $A,$ which connects the sites $x$ and $y$. By definition,
$x\overset{A}{\leftrightarrow}x$ for any $A.$ For every $x\in\Lambda_{n}$ we
define
\[
r_{A}(x)=\sup\{\Vert y\Vert_{\infty}\,:\,y\in\mathbb{Z}^{2}\text{ and
}y\overset{A}{\leftrightarrow}x\}\,.
\]
Clearly, $\Vert x\Vert_{\infty}\leq r_{A}(x)\leq\infty.$ We define
\[
\rho_{V}=\max\{\Vert x\Vert_{\infty};\,x\in V\}\vee1\text{ and }r_{A}
(V)=\max\{r_{A}(x);\,x\in V\}\vee1\,.
\]
Let $R(\delta)$ be the smallest number such that
\[
\mathbb{Q}_{J,\varepsilon}\left(  r_{\mathcal{A}}(V)>R(\delta)\right)
\leq\frac{\delta}{2\Vert f\Vert_{\infty}}\,.
\]
Notice that $R(\delta)<\infty$ since
\[
\mathbb{Q}_{J,\varepsilon}\left(  r_{\mathcal{A}}(V)>R(\delta)\right)
\leq|V|\,\sum_{y:\,\Vert y\Vert_{\infty}>R(\delta)-\rho_{V}}p_{y,\varepsilon
}\,,
\]
and $\sum_{x}p_{x,\varepsilon}=c(\varepsilon)<\infty$, see~(\ref{cepsilon}).

By recurrence of the random-walk $Y_{\cdot}$, which was established above, one
can find, for any $\delta>0$ and $0<\psi<\infty$, a sequence of non-negative
functions $\Psi_{n,\delta,\psi}$ on ${\mathbb{Z}}^{2}$ -- \textit{the
spin-waves} -- such that $\Psi_{n,\delta,\psi}(x)=0$ if $x\not \in\Lambda_{n}
$, $\Psi_{n,\delta,\psi}(x)=\psi$ if $\Vert x\Vert_{\infty}<R(\delta)$, and
\begin{equation}
\lim_{n\rightarrow\infty}\sum_{x\in\Lambda_{n}}\sum_{y\in{\mathbb{Z}}^{2}
}p_{x-y,\varepsilon}\,\left(  \Psi_{n,\delta,\psi}(x)-\Psi_{n,\delta,\psi
}(y)\right)  ^{2}=0\,. \label{conductance}
\end{equation}
The most natural candidate for such a spin-wave is given by
\begin{equation}
\Psi_{n,\delta,\psi}(x)=\psi\,\mathbb{P}_{Y}^{x}\left(  \tau_{\Lambda_{R}
}<\tau_{\Lambda_{n}^{c}}\right)  , \label{voltage}
\end{equation}
where ${\mathbb{P}}_{Y}^{x}$ denotes the law of $Y$-random walk starting at
$x$, whereas $\tau_{\Lambda_{R}}$ and $\tau_{\Lambda_{n}^{c}}$ are the first
hitting times of $\Lambda_{R(\delta)}$ and of the exterior $\Lambda_{n}
^{c}=Z^{2}\setminus\Lambda_{n}$ respectively. Then (\ref{conductance}) is
related to the vanishing, as $n\rightarrow\infty$, of the escape probability
from $\Lambda_{n}$.

The function $\Psi_{n,\delta,\psi}(\cdot)$ in (\ref{voltage}) also represents
the voltage distribution (c.f. \cite{DoS} on the interpretation of recurrence
in terms of electric networks) in the network on the graph $\left(
\mathbb{Z},\mathcal{E}^{J}\right)  $ with bond conductances
$p_{x-y,\varepsilon},$ once all the sites in $\Lambda_{R(\delta)}$ are kept at
the constant voltage $\psi$, whereas all the sites in $\Lambda_{n}^{c}$ are
grounded. In this language the vanishing of the limit in (\ref{conductance})
means zero conductance from $\Lambda_{R(\delta)}$ to infinity, which is a
characteristic property of electric networks corresponding to recurrent random walks.

Let us fix a spin-wave sequence $\{\Psi_{n,\delta,\psi}(x)\}$ so that
(\ref{conductance}) holds.

For any $n$ and\textbf{ }any $A\subset{\mathbb{Z}}^{2}$ such that
$r_{A}(V)\leq R(\delta)$, we define the corresponding $A$-deformed spin-wave
by
\begin{equation}
\widetilde{\Psi}_{n,\delta,\psi,A}(x)\overset{\triangle}{=}\min_{y:x\overset
{A}{\leftrightarrow}y}\Psi_{n,\delta,\psi}(y)\,. \label{84}
\end{equation}
When $A$ is such that $r_{A}(V)>R(\delta)$, we simply set $\widetilde{\Psi
}_{n,\delta,\psi,A}\equiv0$.

For any $x\in\Lambda_{n}$ we denote by $t_{A}(x)\in\mathbb{Z}^{2}$ one of the
sites $y:x\overset{A}{\leftrightarrow}y,$ at which the minimum in (\ref{84})
is attained. (This is a slight abuse of notation, since in fact the site
$t_{A}(x)$ depends also on the function $\Psi_{n,\delta,\psi}(\cdot).$)

The deformed spin-wave is less regular than $\Psi_{n,\delta,\psi},$ but has
the property, crucial for us, that $\widetilde\Psi_{n,\delta,\psi
,A}(x)=\widetilde\Psi_{n,\delta,\psi,A}(y)$ whenever $x\overset{A}
{\leftrightarrow}y.$ In particular, $\widetilde\Psi_{n,\delta,\psi,A}(x)=0$
whenever $x$ is $A $-connected to the outside of $\Lambda_{n}.$

We introduce the tilted measure
\[
\mu_{n}^{U,\bar{\phi},A,\widetilde{\Psi}}(\,\cdot\,)=\mu_{n}^{U,\bar{\phi}
,A}(\,\cdot\,+\widetilde{\Psi}_{n,\delta,\psi,A})\,.
\]
Notice that $\mu_{n}^{U,\bar{\phi},A,\widetilde{\Psi}}=\mu_{n}^{U,\bar{\phi
},A}$ whenever $A$ is such that $r_{A}(V)>R(\delta)$. On the other hand, if
$r_{A}(V)\leq R(\delta)$, then
\[
\left\langle f(\phi+\psi)\right\rangle _{n}^{U,\bar{\phi},A}=\left\langle
f(\phi)\right\rangle _{n}^{U,\bar{\phi},A,\widetilde{\Psi}}\,.
\]
Consequently the following estimate holds:
\begin{align*}
\left|  \left\langle f\left(  \phi+\psi\right)  \right\rangle _{n,\bar{\phi}
}-\left\langle f\left(  \phi\right)  \right\rangle _{n,\bar{\phi}}\right|   &
\leq\mathbb{E}_{n}\left|  \left\langle f(\phi)\right\rangle _{n}^{U,\bar{\phi
},A}-\left\langle f(\phi)\right\rangle _{n}^{U,\bar{\phi},A,\widetilde{\Psi}
}\right| \\
&  +~2\Vert f\Vert_{\infty}\,\mathbb{P}_{n}\left(  r_{A}(V)>R(\delta)\right)
.
\end{align*}
Our target assertion~(\ref{decaylongrange}) is a consequence of the following
two results:
\begin{equation}
\lim_{n\rightarrow\infty}{\mathbb{E}}_{n}\left|  \left\langle f(\phi
)\right\rangle _{n}^{U,\bar{\phi},A}-\left\langle f(\phi)\right\rangle
_{n}^{U,\bar{\phi},A,\widetilde{\Psi}}\right|  =0\,, \label{entropy_bound}
\end{equation}
and
\begin{equation}
2\Vert f\Vert_{\infty}\,\mathbb{P}_{n}\left(  r_{A}(V)>R(\delta)\right)
\leq\delta\,. \label{Pnbound}
\end{equation}
The second bound readily follows from the stochastic domination by the
Bernoulli percolation process $\mathbb{Q}_{J}$ and the definition of
$R(\delta)$. The next subsection is devoted to the proof of \eqref
{entropy_bound}. Our approach is essentially that of~\cite{P,FP}, but with
some simplifications. The main difference between the latter works and ours is
that, using a suitable relative entropy inequality, we obtain estimates on
difference of expectations in finite volume; in this way, \eqref
{entropy_bound} follows immediately by taking the thermodynamic limit, instead
of using the general theory of infinite-volume Gibbs states.

\subsubsection{Relative entropy estimate}

By the well known inequality (see e.g. \cite{F}, f-la (3.4) on p.133),
\[
\left|  \left\langle f(\phi)\right\rangle _{n}^{U,\bar{\phi},A}-\left\langle
f(\phi)\right\rangle _{n}^{U,\bar{\phi},A,\widetilde{\Psi}}\right|  \leq\Vert
f\Vert_{\infty}\sqrt{2\mathrm{H}(\mu_{n}^{U,\bar{\phi},A,\widetilde{\Psi}
}\,|\,\mu_{n}^{U,\bar{\phi},A})},
\]
where $\mathrm{H}(\mu_{n}^{U,\bar{\phi},A,\widetilde{\Psi}}\,|\,\mu
_{n}^{U,\bar{\phi},A})$ is the relative entropy of $\mu_{n}^{U,\bar{\phi
},A,\widetilde{\Psi}}$ with respect to $\mu_{n}^{U,\bar{\phi},A}$. By Jensen's
inequality it suffices to show that
\begin{equation}
\lim_{n\rightarrow\infty}\mathbb{E}_{n}\mathrm{H}\left(  \mu_{n}^{U,\bar{\phi
},{\mathcal{A}},\widetilde{\Psi}}\,|\,\mu_{n}^{U,\bar{\phi},{\mathcal{A}}
}\right)  =0\,. \label{entropy_bound_1}
\end{equation}

From now on we assume that we are working on the event $r_{A}(V)\leq
R(\delta)$ (otherwise the relative entropy is $0$). We follow~\cite{P}, and we
write:
\begin{multline*}
\mathrm{H}(\mu_{n}^{U,\bar{\phi},A,\widetilde{\Psi}}\,|\,\mu_{n}^{U,\bar{\phi
},A})\leq\mathrm{H}(\mu_{n}^{U,\bar{\phi},A,\widetilde{\Psi}}\,|\,\mu
_{n}^{U,\bar{\phi},A})+\mathrm{H}(\mu_{n}^{U,\bar{\phi},A,-\widetilde{\Psi}
}\,|\,\mu_{n}^{U,\bar{\phi},A})\\
=\left\langle \left(  {\mathcal{H}}(\phi+\widetilde{\Psi}_{n,\delta,\psi
,A}\,|\,\bar{\phi})+{\mathcal{H}}(\phi-\widetilde{\Psi}_{n,\delta,\psi
,A}\,|\,\bar{\phi})-2{\mathcal{H}}(\phi\,|\,\bar{\phi})\right)  \right\rangle
_{n,\bar{\phi},A}\,,
\end{multline*}
where, as before, ${\mathcal{H}}(\phi|\bar{\phi})$ is the Hamiltonian defined
by the smooth part of the interaction. Taylor expansion yields
\begin{multline*}
{\mathcal{H}}(\phi+\widetilde{\Psi}_{n,\delta,\psi,A}\,|\,\bar{\phi
})+{\mathcal{H}}(\phi-\widetilde{\Psi}_{n,\delta,\psi,A}\,|\,\bar{\phi
})-2{\mathcal{H}}(\phi\,|\,\bar{\phi})\\
\leq c_{4}\sum_{\substack{x\in\Lambda_{n}\\y\in{\mathbb{Z}}^{2}}
}J_{x-y}\,\left(  \Psi_{n,\delta,\psi}(t_{A}(x))-\Psi_{n,\delta,\psi}
(t_{A}(y))\right)  ^{2}\,
\end{multline*}
with $c_{1}=\max\left|  U^{\prime\prime}\right|  .$ By Jensen's inequality,
\begin{multline}
\left(  \Psi_{n,\delta,\psi}(t_{A}(x))-\Psi_{n,\delta,\psi}(t_{A}(y))\right)
^{2}\\
\leq3\left\{  \left(  \Psi_{n,\delta,\psi}(t_{A}(x))-\Psi_{n,\delta,\psi
}(x)\right)  ^{2}+\left(  \Psi_{n,\delta,\psi}(t_{A}(y))-\Psi_{n,\delta,\psi
}(y)\right)  ^{2}\right. \\
+\left.  \left(  \Psi_{n,\delta,\psi}(x)-\Psi_{n,\delta,\psi}(y)\right)
^{2}\right\}  \,. \label{three}
\end{multline}
The sum of the third terms of (\ref{three}) is bounded by
\[
\sum_{\substack{x\in\Lambda_{n}\\y\in{\mathbb{Z}}^{2}}}J_{x-y}\,\left(
\Psi_{n,\delta,\psi}(x)-\Psi_{n,\delta,\psi}(y)\right)  ^{2}\leq\frac
{1}{\varepsilon}\,\sum_{\substack{x\in\Lambda_{n}\\y\in{\mathbb{Z}}^{2}
}}p_{x-y,\varepsilon}\,\left(  \Psi_{n,\delta,\psi}(x)-\Psi_{n,\delta,\psi
}(y)\right)  ^{2}\,,
\]
and therefore, by the very definition of $\Psi_{n,\delta,\psi}$, goes to zero
as $n\rightarrow\infty$. The contribution of the remaining two terms of
(\ref{three}) to $\mathbb{E}_{n}\mathrm{H}$ is bounded by
\begin{align*}
2\,{\mathbb{E}}_{n}\sum_{\substack{x\in\Lambda_{n}\\y\in{\mathbb{Z}}^{2}
}}J_{x-y}  &  \,\left(  \Psi_{n,\delta,\psi}(t_{A}(x))-\Psi_{n,\delta,\psi
}(x)\right)  ^{2}\\
&  \leq C\,\sum_{x\in\Lambda_{n}}{\mathbb{E}}_{n}\,\left(  \Psi_{n,\delta
,\psi}(t_{A}(x))-\Psi_{n,\delta,\psi}(x)\right)  ^{2}\\
&  \leq C\,\sum_{\substack{x\in\Lambda_{n}\\y\in{\mathbb{Z}}^{2}}
}\mathbb{Q}_{J,\varepsilon}\left(  x\overset{\mathcal{A}}{\leftrightarrow
}y\right)  \left(  \Psi_{n,\delta,\psi}(y)-\Psi_{n,\delta,\psi}(x)\right)
^{2}\\
&  =C\,\sum_{\substack{x\in\Lambda_{n}\\y\in{\mathbb{Z}}^{2}}
}p_{x-y,\varepsilon}\,\left(  \Psi_{n,\delta,\psi}(y)-\Psi_{n,\delta,\psi
}(x)\right)  ^{2}\,,
\end{align*}
and the result follows again from the definition (\ref{conductance}) of
$\Psi_{n,\delta,\psi}$.

\subsection{Continuous symmetry breaking: proof of Theorem \ref{Aiz}}

We construct the states $\mu_{k}$ by prescribing the corresponding boundary
conditions. Let $\Lambda_{n}$ be the box $\left\{  \mathbf{x}\in\mathbb{Z}
^{2}:\left|  \left|  \mathbf{x}\right|  \right|  _{\infty}\leq n\right\}  .$
We define first the boundary condition $\tilde{\phi}_{k}$ by
\[
\tilde{\phi}_{k}\left(  x_{1},x_{2}\right)  =2x_{2}\theta_{k}.
\]
Then it is easy to see that the unique configuration in $\Lambda_{n}$ with
finite energy relative to the b.c. $\tilde{\phi}_{k}$ outside $\Lambda_{n}$ is
the one which coincides with $\tilde{\phi}_{k}$ inside $\Lambda_{n}.$ In
principle that means that the atomic measure concentrated on the configuration
$\tilde{\phi}_{k}$ is itself a Gibbs state for interaction (\ref{85}).
However, this measure is not translation invariant, and also its
finite-dimensional distributions are singular with respect to the Lebesgue
measure. To present a more aesthetically appealing example we proceed as
follows. Denote by $d\nu_{k,\mathbf{x}}^{\delta}$ the uniform distribution on
the circle $\mathbb{S}^{1},$ with the support on the segment of length
$2\delta$ centered around the point $\tilde{\phi}_{k}\left(  \mathbf{x}
\right)  ,$ with $\delta$ small. Denote also by $q\left(  \phi_{\Lambda_{n}
}\Bigm|\mathbf{\phi}\right)  $ the conditional Gibbs distribution in
$\Lambda_{n},$ corresponding to the interaction (\ref{85}) and the boundary
condition $\mathbf{\phi}.$ (In fact, the measures $q\left(  \phi_{\Lambda_{n}
}\Bigm|\mathbf{\phi}\right)  $ are defined only for some b.c. $\mathbf{\phi;}$
namely one needs the set of finite energy configurations to be non-empty. In
the opposite case we define $q\left(  \phi_{\Lambda_{n}}\Bigm|\mathbf{\phi
}\right)  $ to be identically zero measure.) Consider the Gibbs state
$\bar{\mu}_{k,n}^{\delta}$ in $\Lambda_{n},$ given by
\[
\bar{\mu}_{k,n}^{\delta}\left(  d\phi_{\Lambda_{n}}\right)  =Z^{-1}\left(
\Lambda_{n},k,\delta\right)  q\left(  \phi_{\Lambda_{n}}\Bigm|\mathbf{\phi
}\right)  \prod_{\mathbf{x}\in\partial\Lambda_{n}}d\nu_{k,\mathbf{x}}^{\delta
}\left(  \mathbf{\phi}_{\mathbf{x}}\right)  .
\]
Note that for any $\mathbf{x\in}\mathbb{Z}^{2}$ we have that $\bar{\mu}
_{k,n}^{\delta}\left\{  \phi:\left|  \phi\left(  \mathbf{x}\right)
-\tilde{\phi}_{k}\left(  \mathbf{x}\right)  \right|  >\delta\right\}  =0.$

Let $\bar{\mu}_{k}$ be any limit point of the family $\bar{\mu}_{k,n}^{\delta
},$ with $\delta$ small enough, and $\mu_{k}$ be its average over lattice
translations. Then $\mu_{k}$ is a Gibbs state corresponding to the interaction
(\ref{85}), which is not $\mathbb{S}^{1}$-invariant, but which has its
finite-dimensional distributions absolutely continuous with respect to the
Lebesgue measure. $\mathbb{Z}_{k}$-invariance of $\mu_{k}$ is straightforward.

\end{document}